\documentstyle[12pt,twoside]{article}
\def\char{\mbox{char\,}}
\def\Coker{\mbox{Coker\,}}
\def\deg{\mbox{deg\,}}
\def\dim{\mbox{dim\,}}
\def\gd{\mbox{gl. dim\,}}
\def\gkd{\mbox{GK--dim\,}}
\def\gr{\mbox{\bf gr\,}}
\def\id{\mbox{id\,}}

\def\Ker{\mbox{Ker\,}}
\def\len{\mbox{len\,}}

\def\tor{\mbox{Tor\,}}

\def\d{\displaystyle}

\newtheorem{definitia}{Definition}
\newtheorem{prop}{Proposition}
\newtheorem{theorema}[prop]{Theorem}
\newtheorem{lemma}[prop]{Lemma}
\newtheorem{corollarium}[prop]{Corolary}
\def\proof{\noindent {\bf{Proof} }}
\def\nota{\noindent {\bf{Remark} }}

\begin{document}

 \author{D. I. Piontkovsky%
\thanks{
%Moscow State University (Russia 117234, Moskow, MGU, mech--mat faculty,
%High Algebra department)  \newline
       {\sf E-mail: piont@mech.math.msu.su}
                                 }
}

\date{}

\title{Growth and relations in graded rings}

\maketitle

%\noindent

%\bigskip

\section {Introduction}

Let  $k$ be a field. We will call a vector $k$--space, a $k$--algebra,
or
 $k$--algebra  module {\it graded,} if it is ${\bf Z}_+$--graded
and finite--dimensional in every component.
 For a such space  $V$ (in particular, $V$ may be an algebra or a
module)
we denote by $V(x)$ its Hilbert series
$\sum\limits_{i \ge 0} \dim V_i x^i$.
A graded algebra $A= A_0 \oplus A_1 \oplus \dots$
is called {\it connected,} if its zero component $A_0$
is $k$; a connected algebra is called {\it standard,}
if it is generated by $A_1$ and a unit.

The word "algebra" below denotes an  associative  graded algebra.
All inequalities between Hilbert series are coefficient--wise, i.~e.,
we write $\sum_i a_i t^i \ge \sum_i b_i t^i$ iff $a_i \ge b_i$ for all
$i \ge 0$.

%Let $F$ be a free associative algebra and let $I$ be a homogeneous
%ideal in $F$. Known the degrees of generators of $I$, what
%we may say about the growth of an algebra $R = F/I$?
%An answer is given by the Golod--Shafarevitch theorem
%(see~\cite{gsh} or Corollary~\ref{g-sh} below);
%in particular, it follows that if the number of generators of $I$
%having the same degree is sufficiently small, then $R$
%is infinitely dimensional~\cite{gsh} (moreover, if in this case
%$R$ is standard, then it has an exponential growth~\cite{p}).

We are interested in the following situation. Suppose
$A$ is a graded algebra, $\alpha \subset A$ is a subset consisting of
homogeneous elements such that $\alpha$ minimally generates an ideal
$I \triangleleft A$, and $B = A/I$. What can we  say about relations
of the Hilbert series $A(t), B(t)$, and  the generating function
$\alpha (t)$?

If $A$ is a free associative algebra, then a partial answer is given by
Golod--Shafarevich theorem (see~\cite{gsh} or Corollary~\ref{g-sh}
below).
In particular,
it follows that if the number of elements of $\alpha$ of every degree
is sufficiently small, then $B$
is infinite-dimensional~\cite{gsh} (moreover, if in this case
$A$ is standard, then $B$ grows exponentially~\cite{p}).
On the other hand, V.~E.~Govorov (see~\cite{gov} or
Theorem~\ref{govorov}
below) gives some estimates for a number of relations $\# \alpha$,
whenever $B(t)$ is known (if all elements of $\alpha$ have the same
degree).

Our first goal is to demonstrate that the equality cases in
the Golod--Shafarevich theorem and in  both of Govorov's inequalities
 coincide (Theorem~\ref{free}).
Sets $\alpha$ such that  equality holds are called
{\it strongly free}~\cite{an1}, or {\it inert}~\cite{hl};
they are non-commutative analogues of regular sequences
in commutative ring theory \cite{an1}, \cite{golod}.
If $A$ is free, then the set $\alpha \subset A$ of
homogeneous elements is strongly free iff
$$
      \gd A/ id(\alpha) \le 2.
$$
In the general case of an arbitrary connected algebra $A$,
the property of being strongly free is defined by certain conditions on
Hilbert series~\cite{an1}. E.~g., for every set  of
homogeneous elements $\alpha \subset A$ the following inequality holds:
$$
	  (B* k \langle \alpha \rangle) (t) \ge A(t);
$$
the equality holds iff $\alpha$ is strongly free
(where star denotes a free product of algebras).

As for our question in the case of an arbitrary algebra $A$,
we need to consider an asymptotical characteristics of the algebra's
growth.
General graded algebras (even finitely generated ones)
have exponential growth,
so, for a such algebra $A$, $\gkd A = \infty$.
Let us introduce the following analogue of Gelfand--Kirillov dimension:
if $a_i = \dim A_i$, then define the {\it exponent of growth} of $A$ by
$$
      p(A) = \inf \{q>0 \, | \, \exists c>0 \: \forall n\ge 0
		  \quad a_n \le c q^n \}.
$$
At least for a finitely generated algebra $A$, $p(A)$ is finite.
If $A(t)$ is known,  it is clear how to compute $p(A)$:
if $r(A)$ is a radius of convergence of the series $A(t)$, then
$$
    p(A)= \overline{ \lim} \sqrt[n]{a_n} = r(A)^{-1}.
$$
Notice that $p(A)$ depends of the grading on $A$ (for example,
if $A = k \langle x,y | \deg x = \deg y =d \rangle$, then
$p(A)= \sqrt[d]{2}$ depends of $d$), so
it is not a ``dimension'' in the usual sense.

It is proved in~\cite{anick} that over a field of zero characteristic
there is no algorithm to decide whether or not a given quadratic subset
of standard free algebra is strongly free. Using this fact and our
criteria for strongly free sets, we proved in   Section~\ref{rad}
the following. Let $R$ be a finitely presented standard $s$--generated
 algebra with a set of relations $\alpha$.
Then for some particular $s, \alpha$
for some rational numbers $q,r$ there is no algorithm
to decide, when   $\alpha$ is given, whether or not $R(q) = r$.
Moreover, for some integer $n$  there  is no algorithm
to decide, when   $\alpha$ is given, whether or not $p(R) = n$.
It means that in the simplest case that the algebra $A$ is free and
standard,
even the asymptotical version of our question is undecidable in general.

In Section~\ref{extr} below we introduce a concept of {\it extremal
algebra}:
specifically, a graded algebra $A$ is called {\it extremal},
if for any proper
quotient $B = A/I$
$$
	 p(B) < p(A).
$$
If $A$ is an extremal algebra, then it is prime (Theorem~\ref{prime}),
and $p(A)$ is finite (Proposition~\ref{fini}).
Any free product of two connected algebras is extremal
(Theorem~\ref{prod}),
as is any algebra that includes a strongly free set
(Corollary~\ref{sfextr}).

Using extremality, we generalize  Govorov's inequality  for an arbitrary
connected algebra $A$ (Theorem~\ref{main}):
if  $A(t)$ and $B(t)$ are known, we obtain an estimate
for the generating function $\alpha (t)$.
Also, we find a new  characterizing of strongly free sets in terms
of algebras' growth
(Proposition~\ref{r>r}, Theorem~\ref{main}):
if a set $\alpha$ is strongly free, then not only the Hilbert series
$A(t)$, but also the exponent of growth $p(A)$ is as small as possible.
This characterizing generalizes the one of
D.~Anick~\cite[theorem~2.6]{an1}.

I am grateful to Professor E.~S.~Golod for  fruitful discussions.

\section {Golod--Shafarevich theorem and homology}

Let $R = R_0 \oplus R_1 \oplus \dots$ be a standard associative algebra.
Suppose a set $\{x_1, \dots, x_s\}$ is a basis in the space
 $R_1$. Then $R = F/I$, where
$F = k\langle \! x_1, \dots, x_s \! \rangle  = T(R_1)$ is a
free associative algebra, $I \triangleleft F$ is a two--sided ideal
generated by  homogeneous elements of degrees $\ge$ 2.
Let $\alpha = \{f_1, f_2, \dots \}$ be a minimal system of generators of
$I$, and let $u = k \alpha \approx I/ PI+IP$ be a vector space generated
by
$\alpha$
(here
$P = F_1 \oplus F_2 \oplus \dots$ is an augmentation ideal of $F$).

Suppose
  \begin{equation}
  \label{res}
0 \longleftarrow k \stackrel{d_0}{\longleftarrow}
        R \otimes H_0 \stackrel{d_1}{\longleftarrow}
        R \otimes H_1 \stackrel{d_2}{\longleftarrow} \dots
  \end{equation}
be a  minimal free left $R$--resolution  of $k$;
here $H_i \approx \tor_i^R (k,k) = H_i (R),$
i.~e.
$$
  H_0 = k, H_1 \approx R_1, H_2 \approx u
$$
etc.

Let $\Omega^i = \Coker d_{i-1}$ be the $i$--th syzygy module.
Since the Hilbert series of a tensor product is a product
of the factors' Hilbert series,
using the Euler formula for the exact sequence
$$
0 \longleftarrow k \stackrel{d_0}{\longleftarrow}
        R \otimes H_0 \stackrel{d_1}{\longleftarrow}
        \dots \stackrel{d_i}{\longleftarrow} R \otimes H_{i}
	\longleftarrow \Omega^i \longleftarrow 0
$$
we obtain the following
\begin{prop}
\label{resol}
There is an equality of formal power series
$$
 R(x)(1-sx+u(x) - H_3(x)+ \dots + (-1)^i H_i(x)) = 1 + (-1)^i
\Omega^i(x).
$$
In particular, there are inequalities
$$
 R(x)(1-H_0(x)+ \dots - H_{2i-1}(x)) \ge 1,
$$
$$
 R(x)(1-H_0(x)+ \dots + H_{2i}(x)) \le 1;
$$
equalities hold iff
$\gd R \le 2i-1$ (respectively, $\gd R \le 2i$).
\end{prop}

\begin{corollarium}[Golod--Shaferevich theorem]
\label{g-sh}
There is an inequality
  \begin{equation}
  \label{GSH}
    R(x) (1-sx+u(x)) \ge 1.
  \end{equation}
Equality holds iff  $\gd R \le 2$.
\end{corollarium}

The inequality is proved in~\cite{gsh}, the equality condition
is proved in~\cite{an1}. A non-empty set  $\alpha$ such that the
equality
above holds is called  {\it strongly free}, or {\it inert} in $F$;
since in our case $\alpha$ is not empty, these properties
are equivalent to the equality  $\gd R = 2$.

The equality $ R(x) (1-sx+u(x)) = 1 -\Omega^3(x) $ shows that
for a finitely presented algebra $R$ we can
compute the Hilbert series $R(x)$
whenever $\Omega^3(x)$ is known.
So it is interesting to study
the module $\Omega^3$.

\begin{prop}
Denote by $L = F \alpha \triangleleft F$ the left ideal generated by
$\alpha$.
There is an isomorphism of left $R$--modules
$$
   \Omega^3 \cong \tor^F_1 (R, P/L).
$$
\end{prop}

\proof
Suppose $\{u_1, \dots, u_s | \deg u_i = 1\}$ is a basis of the space
$H_1$
and $\{r_1, r_2, \dots | \deg r_i = \deg f_i\}$  is a basis of   $H_2$.
Then we may assume that the map $d_2$ in the resolution~(\ref{res})
is given by the following way:
if $f_i = \sum\limits_{j=1}^n a_i^j x_j,$ then
$d_2 (r_i) = \sum\limits_{j=1}^n \overline {a_i^j} u_j,$
where the overbar denotes the image of an element of  $F$ in $R$.

On the other hand, since any left ideal in a free algebra is a free
module,
taking the long sequence of graded
$\tor_*^F(R, P/L)$ we obtain:
$$
0 \longrightarrow \tor^F_1 (R, P/L) {\longrightarrow}
        R \mathrel{\mathop{\otimes}\limits_F} L
		\stackrel{\phi}{\longrightarrow}
        R \mathrel{\mathop{\otimes}\limits_F} P
		\longrightarrow
        R \mathrel{\mathop{\otimes}\limits_F} P/L
		\longrightarrow  0.
$$
Here
$R \mathrel{\mathop{\otimes}\limits_F} L \cong R \otimes u,$
$R \mathrel{\mathop{\otimes}\limits_F} P \cong R \otimes R_1.$
So  the map $\phi$ induced by the inclusion $L \hookrightarrow P$
coinsides with the map $d_2$ above.
Therefore,
$$
   \tor^F_1 (R, P/L) \cong \Ker d_2 \cong \Omega^3.
$$

\bigskip

%Let us denote $Q = \tor^F_1 (R, P/L)$.

\begin{corollarium}
\label{Q}
The non-empty set $\alpha$ is strongly free iff $\tor^F_1 (R, P/L)=0$.
\end{corollarium}

\nota

All results of this section  hold for an arbitrary connected algebra
$R$.
The only change is to replace in all formulae the term $sx$
by a generating function in the algebra's generators
$\sum_{i \ge 1} t^{\displaystyle \deg x_i}$.

\section {Estimates for the number of relations}

We keep the notation of the previous section.

Suppose that the ideal  $I$ is minimally generated by $t$ ($t>0$)
elements of degree $l$, i.~e.,
$\alpha = \{f_1,  \dots, f_t | \deg f_i = l \}$.
Let  $a_i$ denote the dimension of the space  $R_i$:
$$
       R(x) = \sum_{i \ge 0} a_i x^i.
$$
In the situation above, V.~E.~Govorov proved the following theorem.

\begin{theorema}[\cite{gov}]
\label{govorov}
The series $R(x)$ converges for $x = s^{-1}$.
The following inequalities hold:
\begin{equation}
   \label{Rs}
   t \ge \frac{\displaystyle s^l}{\displaystyle R(s^{-1})}
\end{equation}
and
\begin{equation}
   \label{an}
   t \ge \frac{\displaystyle s a_{n-1} - a_n}{\displaystyle a_{n-l} }
\end{equation}
for all possible  $n$.
Equality holds in~(\ref{an})  for all possible  $n$
if and only if $\tor^F_1 (R, P/L)=0$.
\end{theorema}

The following
theorem shows that the equality hold simultaneously
 in~(\ref{Rs}),~(\ref{an}),
and the Golod--Shafarevich theorem.

\begin{theorema}
\label{free}
The following conditions are equivalent:

  (i)  The set $\alpha$ is strongly free.

  (ii)  Equality holds in~(\ref{Rs}).

  (iii)  Equality holds  in~(\ref{an})   for all possible $n$.
\end{theorema}

\proof

The implications $(i) \Leftrightarrow (iii)$ follow  immediately from
Corollary~\ref{Q} and Theorem~\ref{govorov}.
Let us prove $(i) \Leftrightarrow (ii)$.

Let us take $x=s^{-1}$ in~(\ref{GSH}).
Since $u(x)= t x^l,$  we have the inequality
$$
    R(s^{-1}) t s^{-l} \ge 1,
$$
which is equivalent to~(\ref{Rs}). So, if  $\alpha$ is strongly free,
then equality holds in~(\ref{Rs}).

Conversely, if   $\alpha$  is not strongly free, then the inequality of
formal power series~(\ref{GSH}) is strict, i.~e.,
the following inequality holds
$$
   R(x) (1-sx+tx^l) \ge 1+ a x^n,
$$
where  $n \ge 0, a>0$.
Put  $x = s^{-1}$. We obtain
$$
       R(s^{-1}) t s^{-l} \ge 1+a s^{-n} >1,
$$
or
$$
t > \frac{\displaystyle s^l}{\displaystyle R(s^{-1})}.
$$
So, the inequality~(\ref{Rs}) is strict too.

\section{Radii of convergence of Hilbert series: non-existence of
algorithms}
\label{rad}

                   %Throughout this section, all algebras are connected.
                                        % def. of p(A) and r(A)
For a graded algebra $A$, let $r(A)$  denotes the radius of convergence
of the series $A(t)$.

The following properties of radii of convergence are clear and well
known.

\begin{prop}
\label{r}
Let $A$ be a graded algebra.

(i) $r(A)= \infty$ iff $A$ is finite-dimensional.

(ii) If $A$ is finitely generated, then $r(A) > 0$.

(iii) $r(A) =1$ iff $A$ is not finite-dimensional and $A$ has
      sub-exponential growth.

(iv) If $B$ is either a subalgebra or a quotient algebra of A
with the induced grading, then $r(B) \ge r(A)$.

(v) If $r(A)>0,$ then $\lim \limits_{t \to r(A)} A(t) = \infty$.
So, if $A$ is connected, then the function $f(t) = A(t)^{-1}$
is  continuous  on $[0, r(A)]$  with
$f(0)=1, f(r(A))=0.$
\end{prop}

Under the notation of previous sections, let $l=2$.
It is shown by D.~Anick ~\cite[Theorem~3.1]{anick} that over a field of
zero
characteristic for some positive integers $s$ and $t$, there is no
algorithm
which, when given a set $\alpha \subset F$ of
$t$ homogeneous quadratic elements,
always decides in a finite number of steps whether or not
$\alpha$ is strongly free. We will call such a pair of integers $(s,t)$
{\it undecidable}.

\begin{lemma}
Let $l$ be a positive integer, and let $(s,t)$ be an undecidable pair.
Then the pair $(s+l,t)$ is undecidable.
\end{lemma}

\proof

Let $G = F * k \langle x_{s+1}, \dots, x_{s+l} \rangle$ be
a free algebra of rank $s+l$.
Then a set $\alpha \in F$ is strongly free in $F$ iff it is
strongly free in $G$. So there is no algorithm to recognize it.

\begin{lemma}
Let $s$ be an even integer. If the pair  $(s,t)$ is undecidable,
then the pair $(2s,s^2/4)$ is undecidable too.
\end{lemma}

\proof

By \cite{an4}, a quadratic strongly free set  in $F$
 consting of $q$ elements does exist iff $4q \le s^2$.
So $t \le s^2/4$. Let $\beta$ is a quadratic  strongly free set
in the algebra $G = k \langle x_{s+1}, \dots, x_{2s} \rangle$
consisting of $s^2/4 - t$ elements. Then the set $\alpha \cup \beta$
is strongly free in the algebra
$F*G = k \langle x_{1}, \dots, x_{2s} \rangle$ if and only if
the set $\alpha$ is strongly free in $F$.

\begin{corollarium}
\label{und}
For large enough integer $d$, the pair $(4d,d^2)$ is undecidable.
\end{corollarium}

%Let us denote by

\begin{theorema}

Let $ \mbox{\rm char\,} k = 0$.
Let us denote by  $F_s$ the free associative algebra of rank $s$
with standard grading.
% and let $s,t, F,$ and $R$ be as above.
% Let $A_{s,t}$ be the collection of all quadratic

(i) Let $s,t$ be an undecidable pair of integers.
%For some positive integers$s$ and $t$,
Then there is no algorithm
which, when given a set $\alpha \subset F_s$ of
$t$ homogeneous quadratic elements, always decides whether or not
the equality
$$
R(s^{-1}) = \frac{\displaystyle s^2}{\displaystyle t}
$$
holds, where $R =  F_s/ {\mbox{id\,}}(\alpha)$.

(ii) For some positive integers $s$, $t$, and $q$, there is no algorithm
which, when given a set $\gamma \subset F_{s}$ of
$t$ homogeneous quadratic elements, always decides whether or not
the equality
$$
r (R) = q^{-1}
$$
holds, where $R =  F_{\displaystyle s}/ {\mbox{id\,}}(\gamma)$.
For large enough integer $d$, we can put $s = 64 d$, $t= 241 d^2$, and
$q= 60 d$.
\end{theorema}

\proof

The statement~(i)  follows from Theorem~\ref{free}
and Corollary~\ref{und}.
%and Anick's result above.
To prove~(ii), let $d$ be an integer such that the pair $(60d, 15^2
d^2)$
is undecidable.
Let $\alpha \in F_{60 d}$ be a quadratic set consisting of $15^2 d^2$
elements, and let $B = F_{60 d} / \id (\alpha )$.
Let us denote by $A$ the standard algebra $k\{1,y_1, \dots, y_{4d} \}=
  k \langle y_{1}, \dots, y_{4d} | y_i y_j =0,
     1 \le i,j \le 4d \rangle$.
Then the algebra  $C = k \langle x_1, \dots x_{64d} \rangle / \id
(\gamma)$
is isomorphic to $A*B,$ where
$\gamma = \alpha \cup \{x_i x_j =0 | 60d +1 \le i,j \le 64d\}$.
So $C$ is an algebra with $s=64d$ generators and $t= (4d)^2 + 15^2 d^2 =
241 d^2$ quadratic relations.
It is sufficient to proof that the set $\alpha$ is strongly free in
   $F_{60 d}$ if and only if $r(C) = (60d)^{-1}$.

We have
$$
   C^{-1}(x) = B^{-1}(x) +A^{-1}(x) -1 = B^{-1}(x) + \frac{1}{1 + 4dx}
-1
      = B^{-1}(x) -  \frac{4dx}{1 + 4dx}.
$$
If $\alpha$ is strongly free in
   $F_{60 d}$, then by Corollary~\ref{g-sh} we have $
B^{-1}(x) = 1 - 60d x + 225 d^2 x^2,
$
so
$$
   C^{-1}(x) = 1 - 60d x + 225 d^2 x^2 - \frac{4dx}{1 + 4dx}
%     = \frac{(1 - 60d x)(1 + 4dx) +dx(900 (dx)^2 + 225 dx -4)
%             }{1 + 4dx}$$
%$$
%      = \frac{(1 - 60d x)(1 + 4dx) +dx(60dx -1)(15dx +4)
%             }{1 + 4dx}
      = \frac{(1 - 60d x)(1 - 15 d^2 x^2)
             }{1 + 4dx}.
$$
By Proposition~\ref{r}, $(v)$, we obtain $r(C) = (60d)^{-1}$.

Now let        $\alpha$ is not strongly free in
   $F_{60 d}$.  By Theorem~\ref{free}, we have
$B^{-1}((60d)^{-1}) < 225 d^2 / (60 d)^2 = 1/16$.
%Assume that $r(C) = (60d)^{-1}$; then  $r(B) \ge  (60d)^{-1}$.
%By Corollarium~\ref{g-sh}, we have
%$$ B(x) ( 1 - 60d x + 225 d^2 x^2) > 1, $$
%or
% $$ B(x) ( 1 - 60d x + 225 d^2 x^2) \ge  1+ a x^n, $$
% where $a>0, n>0$.
%If $x$ is a real number such that  $0< x \le (60d)^{-1}$,
%then
%$$   B^{-1}(x) \le \frac{1 - 60d x + 225 d^2 x^2}{1+ a x^n} <
%  1 - 60d x + 225 d^2 x^2, $$
So
$$
  C^{-1}((60d)^{-1}) =
      B^{-1}((60d)^{-1}) -
 \frac{\displaystyle 4d (60d)^{-1}}{\displaystyle 1 + 4d(60d)^{-1}}<
 1/16 -  \frac{\displaystyle 1}{\displaystyle 1+ 1/(4d (60d)^{-1})}= 0.
$$
We obtain  $r(C) \ne (60d)^{-1}$, contradicting  Proposition~\ref{r},
$(v)$.

\section{Algebras of extremal growth}
\label{extr}

Now we introduce the following concept.

\begin{definitia}
A  graded algebra
$A$ is said to be extremal,           % def. of extremal algebras
if for every nonzero homogeneous ideal $I \triangleleft A$
we have $r(A/I) > r(A)$.
\end{definitia}

We will discuss some properties of extremal algebras.

\begin{prop}
\label{fini}
If $A$ is an extremal algebra, then $0 < r(A) < \infty$.
\end{prop}

\proof

%immediately follows from Proposition~\ref{r}, $(i)$.
Let us proof the other inequality.

Suppose that $r(A)=0$.
Let $a \in A$ be a nonzero homogeneous element of degree $d>0$,
let $I$ be an ideal generated by $a$, and let $B=A/I$.
Since $A$ is extremal, $r(B) > 0$.

Let $F=k \langle c|\deg c =d \rangle$.
%By \cite[Lemma~2.2]{an1}, there is an inequality of Hilbert series
By  the obvious inequality of Hilbert series
$$
	  (B*F)(t) \ge A(t),
$$
% (In fact, Anick's proof did not really used an assumption
% that $A$ is connected.)
 $r(B*F)$ must be equal to 0.

On the other hand,
$$
   (B*F)(t)^{-1} = b(t)^{-1} - d t,
$$
where $b(t)$ is equal to $B(t)$ (resp., $B(t)+1$), if   $B$ is unitary
(resp., non-unitary, i.~e. $B_0 = 0$).
Since the right side is an analytical  function in a neighborhood
of zero, and this function takes $0$ into $1$, then in a neighborhood of
zero
its image does not contain $0$.
So the function $(B*F)(t)$ is analitycal in a neighborhood of zero.
Therefore $r(B*F) > 0$.

\begin{theorema}
\label{prime}
Let $A$ be an extremal algebra. Then $A$ is prime.
\end{theorema}

\proof

Obviously, it is sufficient to prove that for any two nonzero
homogeneous
ideals $I$ and $J$ of $A$, $I.J \ne 0$.
Without loss of generality, we can assume that $J$ is a principal
ideal generated by an element $a$ of a degree $h$.

Suppose $I.J = 0$.
Let $B=A/I, C= A/J,$  and let $A = I \oplus V,$ there $V$ is a graded
vector space. We have
$$
   J = k a + A a + a A + AaA = k a + V a + aV +VaV.
$$
Therefore, there is an inequality of Hilbert series
$$
   J(t) \le t^h + 2 t^h V(t) + t^h V(t)^2,
$$
or
$$
   A(t) - C(t) \le t^h (B(t)+1)^2.
$$
So, we have
$$
   A(t) \le C(t) + t^h (B(t)+1)^2.
$$
Thus, for radii of convergence we obtain
$$
   r(A) \ge \min \{r(B), r(C)\},
$$
contradicting the extremality of $A$.

\bigskip

\nota

In fact, we have proved that for any two homogeneous ideals $I,J$
of a graded algebra $A$, if $I.J = 0,$  then
$r(A) = \min \{r(A/I), r(A/J)\}.$

\begin{corollarium}
Let $A$ be a (non-graded) locally finite filtered algebra such that
the associated graded ring $\gr A$ is extremal. Then $A$ is prime.
\end{corollarium}

Now, let us consider examples of extremal algebras.
In fact, the extremality of nontrivial free algebras is proved
by V.~E.~Govorov in~\cite{gov}. We add the following

\begin{theorema}
\label{prod}
Let $A, B$ be non-trivial connected algebras such that $r(A)>0$ and
$r(B)>0$. Then the algebra $A*B$ is extremal.
\end{theorema}

{\bf Proof of Theorem~\ref{prod}}

Let $C=A*B$.
By the formula
$$
 C(t)^{-1} = A(t)^{-1} + B(t)^{-1} -1,
$$
the function $C(t)^{-1}$ is analytical and nonzero in a neighborhood
of zero, so $r(C) > 0$. For $t \in (0, r(C)]$, we have
$0 \le A(t)^{-1} < 1$
and
$0 \le B(t)^{-1} < 1$.
Since $A(r(C))^{-1} + B(r(C))^{-1} -1 = 0,$
we obtain $A(r(C))^{-1}>0$ and $B(r(C))^{-1} >0$;
hence $r(C) < \min \{r(A), r(B) \}$.

It follows from the standard Gr\"obner bases arguments
that we may assume  the algebras $A$ and $B$ to be monomial.
(Indeed, if we fix an order on monomials, then, denoting
by $\overline R$ the associated monomial algebra of an algebra
$R$, we have
$\overline C = \overline A * \overline B$; moreover,
if $I$ is a homogeneous ideal
in $C$, then there exists an ideal $J \triangleleft \overline C$ such
that
$r(C/I) \ge r(\overline C / J)$.)

Now, let $A=k \langle X \rangle /I$ and $B=k \langle Y \rangle /J$,
there $X, Y$ are homogeneous sets
minimally generating algebras $A$ and $B$,
and $I,J$ are ideals generated by monomials of elements of $X$ and $Y$.
If $A$ and $B$ are two--dimensional, i.~e.,
$A \cong k \langle x|x^2=0 \rangle $ and
$B \cong k \langle x|x^2=0 \rangle $, then there is nothing to prove.
So, we can assume that $\dim B \ge 3$.

Let $S \triangleleft C$ be a nonzero principal ideal generated by a
non-empty monomial $m$: it is sufficient to prove
that $r(C/S) > r(C)$.
We will say that a  non-empty  monomial $a$ is an {\it overlap} of two
monomials $b,c$ if there are  non-empty
monomials $f,g$ such that $b=fa, c = a g$.
Now we need the      following

\begin{lemma}
\label{monoms}
Let $A,B$ be connected monomial algebras such that
$\dim_k B \ge 3$, where  $B$ is minimally generated by the set
$Y = \{y_i \}_{\displaystyle i \in \Gamma}$, and let $S$
be a monomial ideal in the algebra $C = A * B$.
Then  $S$ contains
a nonzero monomial $p$ with the following properties:
all monomials $p_{ij}=y_i p y_j, \quad i,j \in \Gamma$
are nonzero, and, moreover, for any two monomials $p_{ij}$ and
$p_{kl}$,  there are no overlaps in the case $j \ne k$
and there is the unique overlap $y_j$ in the case $j=k$.
\end{lemma}

{\bf Proof of Lemma~\ref{monoms}}

It is obvious that $S$ contains a nonzero monomial $n$
such that $n=x n' x$, where $n'$ is a monomial, $x \in X$,
where $X$ is the set of generators of $A$.
To construct such a monomial $p$,
let us consider two cases.

{\it Case 1}
Let $\# Y = 1$, i.~e., $Y=\{y\}$.
Since $\dim B \ge 3$, then $y^2 \ne 0$.
Let $l \ge 0$ be the largest integer satisfying $n = (xy^2)^l n_1$,
where $n_1$ is a monomial.
Put $p = y (xy^2)^q  n_1  (yx)^q $, where $q> \max \{l, \len n_1 +3\}$;
then $p_{11} = y^2 (xy^2)^q  n_1  (yx)^q y$.

Assume that a monomial $a$ is an overlap for the pair $p_{11},p_{11}$.
Then there exist non-empty monomials $c,d$ such that $p_{11} = ca = ad$.
Hence $a$ has the form $(y^2x)^q  f  (xy)^q,$ where $f$ is a monomial;
therefore, $\len c = \len d = \len p_{11} - \len a < q$.
Since $p_{11} = y^2 (xy^2)^q  n_1  (yx)^q y = c (y^2x)^q  f  (xy)^q$,
 $c$ has the form $(y^2 x)^r$ for some $r>0$.
By the maximality of $l$, we have $r=1$, so $\len c = \len d =3$.
On the other hand, it follows from the equality $p_{11} = ad$ that
$d$ has the form $(xy)^l$ for some $l$, so $\len d$ must be even.

{\it Case 2}
Let $\# Y \ge 2$, i.~e., $Y=\{y_1, y_2, \dots \}$.
Put $p = (xy_1)^q  n (y_2 x)^q $, where $q> \len n$.
For some  monomials $p_{ij}$ and $p_{kl}$, suppose $a$ is
an overlap such that
$\len a \ge 2$. Then there exists an overlap of monomials $p,p$.
It means that the set $\{ p \}$ is not {\it combinatorial free},
so, it is not strongly free in a free associative algebra generated by
the
set $X \cup Y$ \cite{an1}. By \cite[Proposition~3.15]{hl},
this means that there exist a non-empty monomial $a$ and a monomial $b$
such that
$$
		  p = aba
$$
(at least in the monomial case, the proof in~\cite{hl} did not really
use the assumption $\char k = 0$).
Then $a$ has the form $a=(xy_1)^q c (y_2 x)^q$ for a monomial $c$,
so $\len a \ge 4 q$ and $\len p \ge 8q,$  contradicting  the choice
of $q$.

\bigskip

Returning to the general proof,
let $P \triangleleft C$ be an ideal generated by all of the monomials
$p_{ij}$, and let $D=C/P$.
Since $P \subset S$,  it is sufficient to prove that $r(D) > r(C)$.
To prove this, we will compute the homology of the algebra $D$ and
obtain
its Hilbert series as the Euler characteristic.

Recall how to compute homologies of a monomial algebra
(see details in~\cite{an2}; we use the terminology of~\cite{ufn}).
Suppose $F$ is a free associative algebra generated by a set $X$,
$I \triangleleft F$ is an ideal minimally generated by a set of
monomials
$U$, and $M$ is the quotient algebra $F/I$.
Let us define a concept of a {\it chain of a rank $n$} and its {\it
tail}.

For $n=0$, every generator $x \in X$ is called a chain of rank $0$;
it coinsides with its tail.
%For $n=1$,
For $n>0,$ a monomial $f = gt$ is called
a chain of rank $n$ and $t$ is called its tail,
if the following conditions hold:
{\it (i)} $g$ is a chain of rank $n-1$;
{\it (ii)} if $r$ is a tail of $g$, then $rt = vu,$
where $v,u$ are monomials and $u \in U$;
{\it (iii)} excluding the word $u$ as the end, there are no subwords
of $rt$ lying in $U$.

Let us denote by
$C_n^M$ the set of chains of rank $n$;
for example,
$C_0^M = X$ and
$C_1^M = U$. Then for all $n \ge 0$ there are  the following
isomorphisms
of graded vector spaces:
$$
              k C_n^M \simeq \tor_{n+1}^M(k,k).
$$
Letting $c_j^i$ denote the number of chains of degree $i$ having a rank
$j$,
consider the generating function
$$
    C^M (s,t) = \sum_{i \ge 0} \sum_{j \ge 0} c_j^i s^j t^i.
$$
Arguing as in Proposition~\ref{resol}, taking the Euler characteristic
of the minimal resolution, we obtain
$$
   M(t)^{-1} = 1 - C^M (-1, t);
$$
the formal power series in the right side does exist since every vector
space $\tor_{n}^M(k,k)$ is concentrated in degrees $\ge n$.

By definition, for all $i \ge 0,$  we have
$$
	    C_i^C = C_i^A \cup C_i^B
$$
and $C_0^D = C_0^C$.
Therefore,
$$
	    C^C (s,t) = C^A (s,t) + C^B (s,t),
$$
so
$$
	    C^D (s,t) = C^A (s,t) + C^B (s,t) + C' (s,t),
$$
there the set  $C'$ consists of  chains that have a subword $p$.
Thus the set $C'_0$ is empty, and $C'_1 = \{ p_{ij} \}$.

Let us prove  that the set $C'$ consists of all monomials
of the form
\begin{equation}
\label{chain'}
	  c_1 p c_2 p \dots p c_n
\end{equation}
for $n \ge 2$, where $c_1, \dots, c_n \in C^B$ .
It is clear that all these monomials are chains of $C'$.
Let us prove the converse.

Since $C_0^B = Y$, this is obvious for chains of rank $1$.
Now, let $f = gt \in C'_n$, where  $g$ is a chain of  lesser rank and
 $t$ is a tail of the chain $f$.
Let $r$ be the tail of $g$.
By induction, we may assume that
$g \in C^B$ or $g$ has the form~(\ref{chain'});
in the second case, $r$ is the tail of $c_n$, or has the form $p y_i,$
where $y_i = c_n \in Y$.
If  $t$ is a word of the alphabet $Y$, then $gt \in C^B$, or
$c_n t \in C^B$, so $f$ has the desired form.
Otherwise, $t$ must contain a subword equal to $p$;
hence, $t = p y_j$ for some $j$.
Thus,
$$
	  f = c_1 p c_2 p \dots p c_n p y_j.
$$

Now, let us compute the generating function.
Notice that if a chain $f$ has the form~(\ref{chain'}),
then the rank of $f$ is equal to $k+n-1,$ where $k$ is the sum of ranks
of the chains $c_1, \dots, c_n$.
Let $\deg p = b$.
By~(\ref{chain'}), we have
$$
   C' (s,t) = \sum_{i \ge 1} (s t^b)^i \left( C^B (s,t) \right)^{i+1} =
    \frac{\d s t^b \left( C^B (s,t) \right)^{2}}{\d 1-s t^b C^B (s,t)}.
$$

Put $q(t) = 1- B(t)^{-1} = C^B (t, -1)$. Obviously,
for $0 < t \le r(B)$,  we have $q(t) > 0$.

We obtain
$$
  D(t)^{-1} = C(t)^{-1} - C'(-1,t),
$$
hence,
$$
  D(r(C))^{-1} = - C'(-1,r(C)) =
          \frac{\d  r(C)^b q(r(C))^2}{\d 1+r(C)^b q(r(C))} > 0.
$$
Since $r(D) \ge r(C)$ and  $D(r(C)) > 0$, we obtain $r(D) > r(C)$.
This completes the proof ot Theorem~\ref{prod}.

\begin{corollarium}
\label{sfextr}
Let $A$ be a connected algebra such that $r(A)>0$.
If there exists a strongly free set in $A$, then $A$ is extremal.
\end{corollarium}

\proof

By \cite[Lemma 2.7]{an1}, any subset of a strongly free set is strongly
free;
so, there is a strongly free element $f \in A$.
Let $L$ be the ideal generated by $f$, and let $B=A/L$.

If $A$ is generated by $f$, then, since every strongly free
set generates a free subalgebra, $A=k \langle f \rangle $;
hence, every proper quotient
of $A$ is finitely-dimensional, so $A$ is extremal.
Otherwise, the algebra $B$ is not trivial,
so the algebra
$C = B * k \langle g| \deg g = \deg f \rangle $ is extremal by
Theorem~\ref{prod}.

By \cite[Section 2]{an1}, there is an isomorphism of graded vector
spaces
$\rho: C \to A$ having the following properties:

(i) the restriction of $\rho$  to $B$ is a right inverse to the
canonical
projection $A \to B$;

(ii) $\rho(g) = f$, and
$ \rho (a_1 g a_2  \dots  g a_n) =
         \rho(a_1) f \rho(a_2)  \dots  f \rho(a_n)$.

Suppose
 $m \in A$ is  an arbitrary homogeneous element,  $m = \rho (c)$, and
$I=AmA$ is the ideal generated by $m$. We need to prove that
$r(A/I) > r(A)$. Indeed, let $c' = g c g$, let $m' = f m f = \rho (c')$,
and let $J \triangleleft A $ (respectively, $K \triangleleft C$)
be the ideal generated by $m'$ (resp., $c'$).
For every $a, b \in C$ we get
$$
	   \rho (a c' b) = \rho (a g c g b) =
               \rho (a) f \rho (c) f \rho (b) =
		    \rho (a) m' \rho (b).
$$
Therefore $\rho (K) \subset J$, so $(A/J)(t) \le (C/K)(t)$.
We obtain
$$
  r(A/I) \ge r(A/J) \ge r(C/K) > r(C) = r(A).
$$

\section{How a quotient algebra may grow?}

Suppose that $A$ is a connected algebra such that $r(A)>0,$
$S \subset A$ is a non-empty set of homogeneous elements minimally
generating an ideal $I=ASA$, and
$B = A/I$.
Let $C = B * k\langle S \rangle$, and let
$$
   D = A/I \oplus I/I^2 \oplus I^2/I^3  \oplus \dots
$$
with the induced grading; then $D(t) = A(t)$.
By~\cite[Theorem~2.4]{hl}, we have a epimorphism
$$
    f: C  \to D,
$$
which is an isomorphism iff the set $S$ is strongly free.
Since $C$
%containes the strongly free subset $S$,
is either a free algebra or a free product of
non-trivial algebras,
it is extremal;
so, we obtain the following

\begin{prop}
\label{r>r}
Using the notation above,
$$
   r(A) \ge r(C).
$$
Equality holds if and only if
the set $S$ is strongly free.
\end{prop}

Now, by the formula for the Hilbert series of a free product,
we have
$$
   C(t)^{-1} = B(t)^{-1} - S(t).
$$
Since $C$ is extremal, the series $B(t)$ and $S(t)$  converge
for $t \in [0, r(C)]$, so
$$
    B(r(C))^{-1} - S(r(C)) =0,
$$
or
$$
    B(r(C)) S(r(C)) =1.
$$
Since $r(A) \ge r(C)$, we have
$$
    B(r(A)) S(r(A)) \ge 1
$$
(where $\infty > 1$); the equality holds iff $r(A) = r(C)$.

Thus we obtain:
\begin{theorema}

\label{main}
Using our notation,
$$
    B(r(A)) S(r(A)) \ge 1,
$$
and the following conditions are equivalent:

(i) the equality above holds;

(ii) $r(A) = r(C)$;

(iii) the set $S \subset A$ is strongly free.
\end{theorema}

In particular, if the set $S$ consists of $t$ elements of degree $l$,
then we have
$$
	 t \ge B(r(A))^{-1} r(A)^{-l},
$$
where  equality holds iff $S$  is strongly free.
This estimate generalizes  Govovrov's inequality~(\ref{Rs}).

\nota

Suppose that $\char k = 0$, the algebra $A$ is free of rank $s$, and
$\alpha$ is a set of $t$ quadratic elements. If the pair $(s,t)$ is
undecidable, then  $C$ is finitely presented and connected
(but non-standard) algebra such that
there is no algorithm do decide  whether or not $r(C) = 1/s$.

\end{document}